\documentstyle[11pt]{article}
\begin{document}

%  Greek letters

\def\a{\alpha}
\def\b{\beta}
\def\d{\delta}
\def\e{\epsilon}
\def\g{\gamma}
\def\k{\kappa}
\def\l{\lambda}
\def\o{\omega}
\def\t{\theta}
\def\s{\sigma}
\def\D{\Delta}
\def\L{\Lambda}

\def\G{{\cal G}}
\def\C{{\bf C}}
\def\P{{\bf P}}

\def\uqoh{U_q(osp(2|2n)^{(1)})}
\def\uqo{U_q(osp(2|2n))}
\def\uqgh{U_q(sl(m|n)^{(1)})}
\def\uqg{U_q(sl(m|n))}
\def\oh{osp(2|2n)^{(1)}}
\def\gh{sl(m|n)^{(1)}}
\def\qh{q^{-h_0/2}}

% Shorthands for \begin{equation} and the like

\def\beq{\begin{equation}}
\def\eeq{\end{equation}}
\def\bea{\begin{eqnarray}}
\def\eea{\end{eqnarray}}
\def\ba{\begin{array}}
\def\ea{\end{array}}
\def\no{\nonumber}
\def\lt{\left}
\def\rt{\right}
\newcommand{\bq}{\begin{quote}}
\newcommand{\eq}{\end{quote}}

\newtheorem{Theorem}{Theorem}
\newtheorem{Definition}{Definition}
\newtheorem{Proposition}{Proposition}
\newtheorem{Lemma}{Lemma}
\newtheorem{Corollary}{Corollary}
\newcommand{\proof}[1]{{\bf Proof. }
        #1\begin{flushright}$\Box$\end{flushright}}

\begin{titlepage}
\begin{flushright}
YITP/K-1113\\
q-alg/9506024
\end{flushright}
\vskip.3in
\begin{center}
{\huge Type-I Quantum Superalgebras, $q$-Supertrace and Two-variable
Link Polynomials}
\vskip.3in
{\Large Mark D. Gould, Jon R. Links}
\vskip.1in
{\large Department of Mathematics, University of Queensland, Brisbane,
Qld 4072, Australia}
\vskip.1in
and
\vskip.1in
{\Large Yao-Zhong Zhang}
\vskip.1in
Yukawa Institute for Theoretical
Physics, Kyoto University, Kyoto 606, Japan

email: yzzhang@yukawa.kyoto-u.ac.jp; yzz@maths.uq.oz.au
\end{center}
\vskip.4in
\begin{center}
{\bf Abstract}
\end{center}
\begin{quote}
A new general eigenvalue formula for the eigenvalues of Casimir invariants,
for the type-I quantum superalgebras, is applied to the construction of
link polynomials associated with {\em any} finite dimensional unitary irrep
for these algebras. This affords a systematic construction of new
two-variable link polynomials asociated with any finite dimensional irrep
(with a real highest weight) for the type-I quantum superalgebras.
In particular infinite families of non-equivalent
two-variable link polynomials are determined in fully explicit form.
\end{quote}

\end{titlepage}
\newpage

\noindent {\bf I. Introduction}
\vskip.1in
Following the celebrated discovery by Jones \cite{Jon85} of the
so-called Jones' link polynomial,
there has been considerable interest in recent years in modern knot theory,
which has  been found to be closely  related, through the
quantum Yang-Baxter equation (QYBE), to various areas of physics such as
solvable models and quantum field theories \cite{Wad89,Wit89}. With
the equally important discovery of quantum algebras during the same period
by Drinfeld \cite{Dri86} and Jimbo \cite{Jim85/86} following the initiatives
of the St. Petersberg group, it was soon
realized by Reshetikhin \cite{Res88} and Turaev \cite{Tur88} that
quantum algebras provided a useful tool in constructing link polynomials.
This idea was further developed in \cite{Zha91,Gou91,Gould} where the authors
proposed a
simple systematic procedure for the construction of link polynomials
arising from  quantum bosonic algebras.

There were many attempts (see e.g. \cite{Wad89,Mur92/93,Aku91})
to construct new two- or multi-variable link
polynomials since the work of HOMFLY \cite{homfly85} and
Kauffman \cite{Kau90} concerning two-variable extensions of the Jones
link polynomial. The two-variable HOMFLY and Kauffman link polynomials
arise from the minimal representations
of $A_n$ and $B_n,~C_n,~D_n$ quantum algebras, respectively.

Subsequently, link polynomials arising from quantum superalgebras
have been addressed
by various authors \cite{Lee90,Deg90,Kau91,Zha92,Lin93,Gou93}.
Among all quantum superalgebras
those of type-I, $U_q[gl(m|n)]$ and $U_q[osp(2|2n)]$,
are particularly interesting because they possess
one-parameter families of finite-dimensional unitary irreps even for
generic $q$. The freedom of having extra parameters in the irreps opens up new
and exciting possibilities in physics \cite{Bra94}.
For the current case, the link
polynomials from such representations will then also depend on these extra
parameters, thus naturally yielding multi-variable link polynomials.
We remark however that such multi-variable link polynomials
are not related \cite{Lin94a}
to those arising from ``colored" braids \cite{Mur92/93}.

For the case of quantum superalgebras, the situation is much more
complicated than the bosonic case. The fundamental difficulty is
the zero $q$-superstrace problem over typical irreps, so that the usual
techniques developed for computing the eigenvalues of Casimir
invariants for quantum bosonic algebras fail in this case.
Due to this problem, only very few isolated examples of multi-variable
link polynomials for quantum superalgebras have so far been known.
These include the two-variable link polynomials
\cite{Lin92} based on $U_q[gl(2|1)]$ and
multi-variable ones \cite{Lin94a}
for a special class of representations of $U_q[gl(m|n)]$. In these
examples, the authors only considered representations for which
the above mentioned difficulty does not occur.

In this paper, we have succeeded in overcoming the above problem
and obtain a well-defined $q$-supertrace formula
which is applied to compute the eigenvalues of the Casimir invariants.
These results are given in theorems \ref{trace-th}, \ref{evalues-th} and
\ref{evalues-th-gen}. However, the proof of the $q$-supertrace formulae
are extremely lengthy and will be published in
a separate paper \cite{Gou95p}. Using these results, we are able to
construct link polynomials associated with
{\em any} finite dimensional unitary irrep of a type-I quantum
superalgebra. Applied to one-parameter families of inequivalent finite
dimensional irreps of $U_q[gl(m|n)]$ and $U_q[osp(2|2n)]$ for generic $q$,
our method affords infinite families of non-equivalent two-variable link
polynomials in fully explicit form.

This paper will be presented in the following order. After recalling
some fundamentals in section II, we give, in section III
an account of the atypicality indices and unitary irreps
of $U_q(\G)$. In section IV we present three theorems
concerning the computation of the $q$-supertraces and
therefore the eigenvalues of Casimir invariants over typical irreps.
Section V derives  a spectral decomposition
formula for the braid generator and its powers.
A general method for constructing link polynomials is presented in
section VI and examples of two-variable link
polynomials are illustrated in section VII.
In the last section, we give a brief discussion of our main results.

\vskip.1in
\noindent {\bf II. Preliminaries}
\vskip.1in
Let $\G$ be a type-I simple Lie superalgebra \cite{Kac77/78}
with generators
$\{e_i,~f_i,~h_i\}$ and let $\a_i,~i=0, 1,\cdots, r,$ be its simple roots
with $\a_0$ the unique odd simple root; here we choose the distinguished
\footnote{Superalgebras allow many
inequivalent systems of simple roots. See \cite{Kac77/78}. The relation
between the different quantum superalgebras obtained by choosing
different systems of simple roots is studied in \cite{Kho94}.}
set of simple roots.
Let $(~,~)$ be a fixed invariant bilinear form on $H^*$, the dual of the Cartan
subalgebra $H$ of $\G$. The quantum superalgebra $U_q(\G)$ has the
structure of a ${\bf Z}_2$-graded quasi-triangular Hopf algebra. Throughout the
paper we will assume that $q$ is generic, i.e. not a root of unity.
We will not give the full defining relations of $U_q(\G)$ here but mention
that the simple raising and lowering generators of $U_q(\G)$ obey more
relations than just the usual $q$-Serre relations known from quantum
bosonic algebras \cite{Kho91,Yam91,Sch92/93,Del94a}.
These necessary extra relations are  referred to
as ``extra $q$-Serre relations". $U_q(\G)$ has a coproduct $\D$ and
antipode $S$ given by
\begin{eqnarray}
&&\Delta(q^{\pm h_i})=q^{\pm h_i}\otimes q^{\pm h_i},\nonumber\\
&&\Delta(e_i)=e_i\otimes q^{-h_i/2}+q^{h_i/2}\otimes e_i,\nonumber\\
&&\Delta(f_i)=f_i\otimes q^{-h_i/2}+q^{h_i/2}\otimes f_i,\label{coproduct}\\
&&S(a)=q^{-h_\rho}\g(a)q^{h_\rho}\,,~~~~a=e_i,~f_i,~h_i\label{antipode}
\end{eqnarray}
where $\g$ is the principal anti-automorphism on $U_q(\G)$ and
$\rho$ is the graded half-sum of positive roots of $\G$.
We omit the formulas for the counit which are not needed here.

%Let $\dagger$ denotes the usual conjugation operation defined on
%$U_q(\G)$, which satisfies
%\beq
%\g(a)^\dagger=\g(a^\dagger),~~~~\forall a\in U_q(\G)
%\eeq
%which implies that $S(a^\dagger)=S^{-1}(a)^\dagger$. Moreover,
%$\D(a)^\dagger=\D(a^\dagger)$, where $\dagger$ on the left hand side is the
%usual conjugation operation induced on $U_q(\G)\otimes U_q(\G)$, so that
%for homogenous elements $a,~b\in U_q(\G)$
%\beq
%(a\otimes b)^\dagger=(-1)^{[a][b]}\;a^\dagger\otimes b^\dagger
%\eeq

The algebra $U_q(\G)$ is a quasitriangular graded Hopf algebra,
which means the following.
Let $\D'$ be the opposite coproduct so that
$\D'=T\D$, where $T$ is the graded twist
map: $T(a\otimes b)=(-1)^{[a][b]}b\otimes a\,,~\forall a,b\in U_q(\G)$.
Here $[a]\in{\bf Z_2}$ denotes the grading of element $a:~[a]=0$ if $a$
is even and $[a]=1$ if it is odd.
Then $\Delta$ and $\Delta'$ are related by the universal R-matrix $R$
in $U_q(\G)\otimes U_q(\G)$ satisfying, among others,
the relations
\begin{eqnarray}
&&R\D(a)=\D'(a)R\,,~~~~~\forall a\in U_q(\G)\label{intw}\\
&&(I\otimes \D)R=R_{13}R_{12}\,,~~~~(\D\otimes I)R=R_{13}R_{23}\label{hopf}
\end{eqnarray}
where if $R=\sum a_t\otimes b_t$ then $R_{12}=\sum a_t\otimes b_t\otimes 1$
, $R_{13}=\sum a_t\otimes 1\otimes b_t$ etc. It follows from (\ref{hopf})
that $R$ satisfies the QYBE:
\beq
R_{12}R_{13}R_{23}=R_{23}R_{13}R_{12}.\label{uqybe}
\eeq
Note that the multiplication rule for the tensor product is defined for
homogeneous elements $a,b,c,d\in U_q(\G)$ by
\begin{equation}\label{gradprod}
(a\otimes b)(c\otimes d)=(-1)^{[b][c]}(ac\otimes bd).
\end{equation}

It is a well established fact for
quasitriangular Hopf algebras, that there exists a distinguished element
\cite{Dri86}
\begin{equation}
u=\sum_t(-1)^{[t]}S(b_t)a_t\label{u}
\end{equation}
where, as above,  $a_t$ and $b_t$ are coordinates of the universal $R$-matrix
. One can show that $u$ has inverse
\begin{equation}
u^{-1}=\sum_t(-1)^{[t]}S^{-2}(b_t)a_t\label{u-1}
\end{equation}
and satisfies
\begin{eqnarray}
&&S^2(a)=uau^{-1}\,,~~~\forall a\in U_q(\G)\nonumber\\
&&\Delta(u)=(u\otimes u)(R^TR)^{-1}
\end{eqnarray}
where $R^T=T(R)$.  It is easy to check that
\beq
v=u^{-1}q^{-2h_\rho}
\eeq
belongs to the center of $U_q({\cal G})$ and satisfies
\begin{equation}
\Delta(v)=(v\otimes v) (R^TR)^{-1}\label{vv1}.
\end{equation}
Moreover, on a finite dimensional irreducible module $V(\L)$
with highest weight
$\L\in D^+$, the Casimir operator $v$ takes the eigenvalue
\begin{equation}
\chi_\L(v)=q^{(\L,\L+2\rho)}.\label{chi1}
\end{equation}

Note that the generators $\{e_i,~f_i,~q^{h_i},~i=1,\cdots,r\}$ form generators
of the quantum group $U_q(\G_0)$, where $\G_0$ is the
``even subalgebra" of $\G$. Specifically,
\bea
&&\G_0=u(1)\bigoplus sl(m)\bigoplus sl(n)\,,~~~~{\rm for}~\G=sl(m|n),
{}~~~m,n\geq2,\no\\
&&\G_0=u(1)\bigoplus sl(n)\,,~~~~{\rm for}~\G=sl(1|n),
{}~~~n\geq2,\no\\
&&\G_0=u(1)\bigoplus sp(2n)\,,~~~~{\rm for}~\G=osp(2|2n).
\eea
Throughtout we let  $V_0(\L)$ denote the finite dimensional irreducible
$U_q(\G_0)$ module with higest weight $\L\in D^+$. We call
\begin{equation}
D^0_q(\Lambda)=\prod_{\b\in\Phi^+_0}\frac{[(\Lambda+\rho,\b)]_q}
  {[(\rho,\b)]_q}\label{q-dimension}
\end{equation}
the $q$-dimension of the $U_q(\G_0)$ irrep $V_0(\L)$,
where $\Phi^+_0$ denotes the set of even positive
roots of ${\cal G}$. Here and in what follows we will adopt the notation
\beq
[x]_q=\frac{q^x-q^{-x}}{q-q^{-1}}
\eeq

\vskip.1in
\noindent{\bf III. Atypicality Indices and Finite-Dimensional Unitary Irreps}
\vskip.1in
Let $K(\L)$ be the Kac-module associated to $V(\L)$. $K(\L)$
is not necessarily irreducible. If
it is, we have $V(\L)=K(\L)$
and refer to $\L$ and $V(\L)$ as ``typical".
Recall that $\L$ is typical iff $(\L+\rho,\b)\neq 0,~~\forall \b
\in \Phi_1^+$, where
$\Phi_1^+$ is the set of odd positive roots of $\G$.

Let us remark that for typical modules the dimensions are easily evaluated
to be ${\rm dim}V(\L)=2^d\cdot{\rm dim}V_0(\L)$,
where $d$, which is equal to $mn$ for $gl(m|n)$ and $2n$ for $osp(2|2n)$,
is the number of odd positive roots. This formula is
particularly useful in determining tensor product decompositions of typical
modules.

\begin{Definition}\label{atypical-def}
The integer
\beq
a_\L=\left |\bar{\Phi}_1^+(\L)\right |,~~~~\bar{\Phi}_1^+(\L)
  =\{\b\in\Phi_1^+|(\L+\rho,\b)=0\}
\eeq
is called the ``atypicality index" of $\L\in D^+$.
In particular, $a_\L=0$ iff $\L$ is typical.
\end{Definition}

The type-I quantum superalgebras admit two types of unitary representations
which may be described as follows. We make the simplifying assumption
that $q>0$ (i.e. $q$ is real and positive) and define a conjugation operation
on the $U_q(\G)$ generators by
$e_i^\dagger=f_i\,,~~~~f_i^\dagger=e_i\,,~~~~h_i^\dagger=h_i$
which is extended uniquely to all of $U_q(\G)$ such that
$(xy)^\dagger=y^\dagger x^\dagger\,,~~~~~~\forall x,y\in U_q(\G)$.
We call $\pi_\L$ type (1) unitary if
\beq
\pi_\L(x^\dagger)=\overline{\pi_\L(x)}\,,~~~~\forall x\in U_q(\G)
\eeq
and type (2) unitary if
\beq
\pi_\L(x^\dagger)=(-1)^{[x]}\;\overline{\pi_\L(x)}\,,~~~~
  \forall x\in U_q(\G)
\eeq
where the overline denotes Hermitian matrix conjugation.
The two types of unitary representations are in fact related via duality.

\begin{Lemma}\label{lemma-unitary}
Such unitary representations have the property that they are always
completely reducible and the tensor product of
two irreducible unitary representations of the
same type reduces completely into irreducible
unitary representations of the same type.
Moreover the atypicality indices of the irreps occurring in this decomposition
are less than or equal to the atypicality index of either component.
\end{Lemma}

The finite dimensional irreducible unitary representations for all
type-I quantum superalgebras have been
classified in \cite{Gou95,Lin95}. For completeness we cite these
classification results below. Let us first
of all introduce some notation. For $gl(m|n)$,  we choose
$\{\e_i\}^m_{i=1}\bigcup \{\d_j\}^n_{j=1}$ as a
basis for $H^*$ with $[\e_i]=0,~[\d_j]=1$ and
\begin{equation}
(\e_i,\e_j)=\delta_{ij},~~~~
(\d_i,\d_j)=-\delta_{ij},~~~~ (\e_i, \d_j)=0\,.
\end{equation}
Using this basis, any weight $\Lambda$ may written as
\begin{equation}
\Lambda\equiv (\Lambda_1,\cdots,\Lambda_m|\bar{\Lambda}_1,
\cdots, \bar{\Lambda}_n)\equiv \sum_{i=1}^m\Lambda_i\e_i+
\sum_{j=1}^n\bar{\Lambda}_j\d_j
\end{equation}
and the graded half-sum $\rho$ of the positive roots is
\begin{equation}
2\rho=\sum_{i=1}^m(m-n-2i+1)\varepsilon_i+\sum_{j=1}^n(m+n-2j+1)
\d_j\,.
\end{equation}
For $osp(2|2n)$,
choose $\{\e_0\}\bigcup\{\e_i\}_{i=1}^n$ as a basis for $H^*$
with $[\e_0]=1,~[\e_i]=0$ and
\beq
(\e_0,\e_0)=-1\,,~~~(\e_i,\e_j)=\d_{ij}\,,~~~\forall i, j=1,\cdots, n\,,~~~
  (\e_0,\e_i)=0.
\eeq
In this case, any weight $\L$ may be expressed as
\begin{equation}
\Lambda\equiv (\bar{\L}|\Lambda_1,\cdots,\Lambda_n)\equiv \bar{\L}\e_0+
\sum_{i=1}^n\Lambda_i\e_i
\end{equation}
and the graded half-sum $\rho$ of the positive roots is given by
\beq
\rho=\sum^n_{i=1}(n-i+1)\e_i-n\e_0.
\eeq

\begin{Proposition}\label{prop1}
(I) A given $U_q[gl(m|n)]$-module
$V(\L)$, with $\L\in D_+$, is type (1) unitary iff: i) $(\L+\rho,\e_m
-\d_n)>0$; or ii) there exists an odd index $\omega\in\{1,2,\cdots,n\}$
such that $(\L+\rho,\e_m-\d_\omega)=0=(\L,\d_\omega-\d_n)$.
In the former case the given condition also enforces typicality on
$V(\L)$, while in the latter case all irreps are atypical.\\
(II) The $U_q[gl(m|n)]$-module $V(\L)$, with $\L\in D_+$, is type (2)
unitary iff: i) $(\L+\rho,\e_1
-\d_1)<0$; or ii) there exists an even index $k\in\{1,2,\cdots,m\}$
such that $(\L+\rho,\e_k-\d_1)=0=(\L,\e_1-\e_k)$.
In the former case $V(\L)$ is typical, while in the latter case it is atypical.
\end{Proposition}

\begin{Proposition}\label{prop3}
A given $U_q[osp(2|2n])$-module $V(\L)$ is type (1) unitary
iff $(\L,\a_0)\geq 0$,
where $\a_0$ denotes the unique odd simple root, and type (2) unitary iff
i) $(\L+\rho, \e_0+\e_1)<0$; or ii) there exists an index
$k\in \{1,2,\cdots,n\}$
such that $(\L+\rho,\e_0+\e_k)=0=(\L,\e_1-\e_k)$; or iii) $\L=0$.
\end{Proposition}

\vskip.1in
\noindent {\bf IV. $q$-Supertrace and Eigenvalues of Casimir Invariants}
\vskip.1in
Throughtout this section we assume $V(\L)$ is a fixed but arbitrary
finite dimensional irreducible $U_q(\G)$-module.
Suppose $V(\nu)\subset V(\mu)\otimes V(\L)$ is typical and
let $P[\nu]$ denote the central projection of the
tensor product module $V(\mu)\otimes V(\L)$ onto its isotypic component
$\bar{V}(\nu)\equiv m_\nu V(\nu)\subset V(\mu)\otimes V(\L)$
(that is $\bar{V}(\nu)=V(\nu)\oplus\cdots\oplus V(\nu),~~
m_\nu$ copies).
We state the following $q$-supertrace formula:

\begin{Theorem}\label{trace-th}
For $\mu,~\nu\in D^+$ typical,
\beq
(I\otimes {\rm str})(I\otimes\pi_\L(q^{-2h_\rho}))P[\nu]=(-1)^{[\nu]}\,
  m_\nu\,\frac{\chi_\mu(\Gamma_0)}{\chi_\nu(\Gamma_0)}\cdot
  \frac{D^0_q(\nu)}{D^0_q(\mu)}
\eeq
where $[\nu]$ modulo 2 is the degree of the weight $\nu$,
$\Gamma_0$ is a central element of $U_q(\G_0)$ and $\chi_\mu(\Gamma_0)$
is the eigenvalue of $\Gamma_0$ on the $U_q(\G_0)$-module $V_0(\mu)$:
\beq
\chi_\mu(\Gamma_0)=
   \prod_{\b\in\Phi^+_1}\frac{[(\mu+\rho,\b)]_q}
   {[(\rho,\b)]_q}.\label{chi-gamma}
\eeq
\end{Theorem}

The proof of this theorem is very lengthy and detailed,
and will be published elsewhere \cite{Gou95p}.

\begin{Proposition}
If the operator $c\in U_q(\G)\otimes {\rm End}V(\L)$ satisfies
$\D_\L(a)c=c\D_\L(a),~~\forall a\in U_q(\G)$, where $\D_\L=(I\otimes
\pi_\L)\D$, then
\begin{equation}
C^\L_k=(I\otimes {\rm str})\{[I\otimes \pi_\Lambda(q^{-2h_\rho})]c^k\}
\,,~~~~~k\in{\bf Z}^+\label{casimir-gen}
\end{equation}
belong to the center of $U_q(\G)$ and thus form a family of Casimir invariants
{}.
\end{Proposition}

An important example of $c$ is given by
\beq
c=\frac{I\otimes I-R^T_\L R_\L}{q-q^{-1}}
\eeq
where $R_\L=(I\otimes \pi_\L)R$, with $R$ the universal R-matrix.

Now assume $V(\mu)\otimes V(\L)$ is completely reducible and write
\beq
V(\mu)\otimes V(\L)=\bigoplus_{\nu}m_\nu V(\nu)\label{decomposition-gen}
\eeq
with now $m_\nu$ the multiplicity of the module $V(\nu)$ occurring in the
tensor
product. This always occurs when $\mu$ and $\L$ are unitary of the
same type. Moreover, in such a case, each of the modules $V(\nu)$ is also
unitary. If $c\in (I\otimes \pi_\L)(Z\otimes Z)\D(Z)$
where $Z$ is the centre of $U_q(\G)$, then one
can deduce the following spectral decomposition for
$c$ and its powers $c^k,~k\in{\bf Z}$:
\beq
c^k=\sum_\nu\left (\chi_\nu(c)\right )^k P[\nu]
\eeq
where $\chi_\nu(c)$ is the eigenvalue of $c$ on $V(\nu)\subset V(\mu)\otimes
V(\L)$. Thus if $c$ is given by the above example, then we have
\beq
\chi_{\nu}(c)=\frac{1-q^{C(\mu)+C(\L)-C(\nu)}}{q-q^{-1}}
\eeq
where $C(\L)\equiv (\L,\L+2\rho)$ denotes the eigenvalue of the second order
Casimir invariant of $\G$.

With the aid of theorem \ref{trace-th}, we can determine the eigenvalues of
the Casimir invariants $C^\L_k$ on a finite dimensional typical module
$V(\mu)$ [notation as in (\ref{decomposition-gen})].

\begin{Theorem}\label{evalues-th}
If $\mu,~\nu$ are all typical, then the eigenvalues of the Casimir invariants
on $V(\mu)$ are given by
\begin{equation}
\chi_\mu(C^\L_k)=\sum_{\nu}(-1)^{[\nu]}\,m_\nu\,
  \left(\chi_{\nu}(c)
  \right)^k \frac{\chi_\mu(\Gamma_0)}{\chi_{\nu}(\Gamma_0)}
\cdot \frac{D^0_q(\nu)}{D^0_q(\mu)}
\,,~~~~k\in {\bf Z}.\label{chi-casimir}
\end{equation}
\end{Theorem}

\vskip.1in
\noindent{\bf Remark:}
Let  $\{\l_i\}$  denote the set of
distinct weights in $V(\L)$ occurring with multiplicites $m_{\l_i}$.
It can be shown that the
above theorem may be extended to all finite dimensional modules
$V(\mu),~\mu\in D^+$, by replacing $\nu,~m_\nu$
with $\mu+\l_i,~m_{\l_i}$ respectively,  and summing
over $\l_i$. For more details see \cite{Gou95b}. The eigenvalue
formula obtained in this way is referred to as the ``extended eigenvalue
formula" on $V(\mu),~\mu\in D^+$.
Note that for generic $\mu$, the extended eigenvalue formula
determines a polynomial function on
$H^*$. It is well defined if all $\mu+\l_i$ are typical but if some
$\mu+\l_i$ is atypical it is necessary first to expand the right hand side
of the extended eigenvalue formula into a polynomial
in order to avoid singularities \cite{Gou95b}.

In the case of
unitary $\mu\in D^+$ this latter problem can be overcome as follows. We set
\beq
\Phi_1^+(\l)=\{\b\in\Phi_1^+|(\l+\rho,\b)\neq 0\}
\eeq
so that
\beq
\left |\Phi_1^+(\l)\right | +a_\l=\left |\Phi_1^+\right |.
\eeq
Then we have the following

\begin{Theorem}\label{evalues-th-gen}
The eigenvalues of the Casimir invariants on a unitary module
$V(\mu)$ are given by
\begin{equation}
\chi_\mu(C^\L_k)=\sum_{\{\nu\,|\,a_{\nu}=a_\mu\}}(-1)^{[\nu]}\,m_\nu\,
  \left(\chi_{\nu}(c)
  \right)^k \frac{\prod_{\b\in\Phi_1^+(\mu)}\;[(\mu+\rho,\b)]_q}
  {\prod_{\b\in\Phi_1^+(\nu)}\;[(\nu+\rho, \b)]_q}
\cdot \frac{D^0_q(\nu)}{D^0_q(\mu)}
\,,~~~k\in{\bf Z}^+\label{chi-casimir-gen}
\end{equation}
provided that $V(\L), ~V(\mu)$ are unitary of the same type. Here the
sum over $\nu$ is over $V_0(\nu)\subset V_0(\mu)\otimes V(\L)$ and $m_\nu$
is the multiplicity of $V_0(\nu)$ in this space.
\end{Theorem}

For a given unitary module $V(\L)$,
the above formula is well defined for all unitary $\mu\in D^+$ of the
same type.

\vskip.1in
\noindent {\bf V. Diagonalization of the Braid Generator}
\vskip.1in
Let $P$ be the graded permutation operator on $V(\L)\otimes V(\L)$ defined by
$P(|x>\otimes |y>)=(-1)^{[x][y]}
|y>\otimes |x>$, for all homogeneous $|x>\,,\,|y>\in V(\L)$ and set
\begin{equation}
\sigma=PR\,~~~~~\in\,{\rm End}(V(\L)\otimes V(\L)).\label{pr}
\end{equation}
Here and in what follows we regard elements of $U_q(\G)$ as
operators on $V(\L)$. Then (\ref{intw}) is equivalent to
\begin{equation}
\sigma\Delta(a)=\Delta(a)\sigma\,~~~~\forall a\in U_q(\G)\label{sigma-delta}
\end{equation}
and (\ref{uqybe}) can be written as
\beq
(I\otimes \s)(\s\otimes I)(I\otimes\s)=(\s\otimes I)(I\otimes\s)
   (\s\otimes I).\label{ybe-braid}
\eeq
It follows immediately that the operaotrs $\s_i^\pm\in {\rm End}
(V(\L)^{\otimes M})\,~~i=\{1,2,\cdots, M-1\}$ defined by
\beq
\s_i^\pm=\underbrace{I\otimes\cdots\otimes I}_{i-1}
  \otimes \s^\pm\otimes\underbrace{I\otimes\cdots\otimes I}_{M-i-1})
  \label{braid}
\eeq
generate a nontrivial representation of the rank $(M-1)$ braid group $B_M$.

In the case when $\s$ acts on $V(\L)\otimes V(\L)$ with $V(\L)$ unitary,
it can be shown that $\s$ is self-adjoint
and diagonalizable \cite{Lin94b}. We remark
however that only the type-I quantum superalgebras admit finite
dimensional unitary irreps.

Similar to (\ref{decomposition-gen}) we write,
\beq
V(\L)\otimes V(\L)=\bigoplus_{\nu}m_\nu V(\nu)\label{decomposition}
\eeq
where again $m_\nu$ is the multiplicity of
the module $V(\nu)$ occurring in the tensor
product and each of the modules $V(\nu)$ is unitary. In
view of the self-adjointness of $\s$, $\s$ is diagonalizable on
$\bar{V}(\nu)\equiv m_\nu V(\nu)=V(\nu)\bigoplus\cdots\bigoplus V(\nu)~~
  (m_\nu ~{\rm copies})$,
regardless of the multiplicity. In fact it is possible to derive a
spectral decomposition formula for $\s$, as in the case of quantum bosonic
algebras \cite{Gould}.

Recall that $\lim_{q\rightarrow 1}\;\sigma=P$ and $P$ is diagonalizable
on $V(\L)\otimes V(\L)$ with eigenvalues $\pm 1$. Following \cite{Gould},
let $P[\pm]$ denote the projection operators defined by
\begin{equation}
P[\pm](V(\L)\otimes V(\L))=W_\pm
\end{equation}
where
\begin{equation}
W_\pm=\{w\in V(\L)\otimes V(\L)\,|\, \lim_{q\rightarrow 1}(\sigma\mp 1)w=0\}.
\end{equation}
Since $\sigma$ is an $U_q(\G)$-invariant each subspace $W_\pm$
determines a $U_q({\cal G})$-module and $P[\pm]$ commute with the action
of $U_q({\cal G})$.
As above $P[\nu]$ denotes the projection operator onto the modules
$\bar{V}(\nu)$; then obviously
\beq
P[\nu,\pm]=P[\pm]P[\nu]=P[\nu]P[\pm]
\eeq
is the projection onto the isotypic component $\bar{V}(\nu)$
consisting of eigenvectors of $\s$ with parities $\pm 1$, respectively
(i.e. the component of $\bar{V}(\nu)$ in $W_\pm$ respectively).

The diagonalizability of $\s$, together with the fact that
\begin{equation}
\sigma^2=PRP\cdot R=R^TR=(v\otimes v)\Delta(v^{-1}),
\end{equation}
implies the following spectral decomposition for $\s$ and its powers:
\beq
\s^k=q^{-kC(\L)}\sum_{\nu} q^{\frac{k}{2}C(\nu)}\left (P[\nu,+]+
  (-1)^k P[\nu,-]\right ),~~~~k\in {\bf Z}\label{spectral}
\eeq
where as before $C(\l)\equiv (\l, \l+2\rho)$.
It follows in particular that $\s$ satisfies the polynomial identity
\beq
\prod_\nu\left (\s-q^{\frac{1}{2}C(\nu)-C(\L)}\right )
   \left (\s+q^{\frac{1}{2}C(\nu)-C(\L)}\right )=0
\eeq
which leads to the generalized skein relations for the corresponding link
polynomials investigated below.

\vskip.1in
\noindent {\bf VI. Link Polynomials}
\vskip.1in
Let $\t\in B_M$ be a word in the generators $\s_i^\pm,~~1\leq i\leq M-1$
and let $\hat{\t}$ denote the link obtained by closing the braid.
For the construction of link polynomials, the Markov trace $\phi$
plays an essential role. It is defined by
\bea
(i)~~~~~~~&&\phi(\t\eta)=\phi(\eta\t),~~~~~~\forall \t,
   \eta\in B_M\nonumber\\
(ii)~~~~~~&&\phi(\t\s_{M-1})=z\phi(\t),~~~~\phi(\t\s_{M-1}^{-1})=\bar{z}
   \phi(\t),~~~~\forall \t\in B_{M-1}\subset B_M.
\eea
Given such a Markov trace, it is well-known that one can define a link
polynomial $L(\hat{\t})$ through
\beq
L(\hat{\t})=(z\bar{z})^{(M-1)/2}\;(\bar{z}z^{-1})^{e(\t)/2}\phi(\t),~~~~
   \t\in B_M
\eeq
where $e(\t)$ is the sum of the exponents of the $\s_i$'s appearing in $\t$.
The functional $L(\hat{\t})$ enjoys the following properties:
\bea
(i)~~~~~~~&&L(\widehat{\t\eta})=L(\widehat{\eta\t}),~~~~~~\forall \t,
   \eta\in B_M\nonumber\\
(ii)~~~~~~&&L(\widehat{\t\s_{M-1}^{\pm 1}})=L(\hat{\t}),~~~~~~
   \forall \t\in B_{M-1}\subset B_M
\eea
and is an invariant of ambient isotopy.

\begin{Proposition}\label{link-prop}
The functional $\phi(\t)$ defined by
\beq
\phi(\t)=\frac{({\rm tr}\otimes {\rm str}^{\otimes (M-1)})(I\otimes
  \D^{(M-1)}(q^{-2h_\rho})\t)}{{\rm dim} V(\L)}
\eeq
where tr and str denote the trace and supertrace over $V(\L)$ respectively,
qualifies as a Markov trace with
\beq
z=q^{(\L,\L+2\rho)},~~~~~\bar{z}=q^{-(\L,\L+2\rho)}.
\eeq
\end{Proposition}

\begin{Corollary}\label{link-corollary}
It follows that
\beq
L(\hat{\t})=q^{-(\L,\L+2\rho)e(\t)}\phi(\t),~~~~~~~\t\in B_M
\eeq
defines a link polynomial.
\end{Corollary}

Now consider the family of Casimir invariants
\beq
C^\L_k=(I\otimes {\rm str})[I\otimes \pi_\L(q^{-2h_\rho})]\s^k.
\eeq
Let $\xi_k^\L$ denote the eigenvalues of the invariants $C_k^\L$ on
$V(\L)$. In view of (\ref{spectral}) and theorem \ref{evalues-th},
one can deduce, for $\L$ typical, that they are given explicitly by
\beq
\xi^\L_k=q^{-kC(\L)}\sum_\nu(-1)^{[\nu]}\,q^{\frac{k}{2}
  C(\nu)}\left (m_\nu^++(-1)^km_\nu^-\right )
  \frac{\chi_\L(\Gamma_0)}{\chi_{\nu}(\Gamma_0)}
  \cdot \frac{D^0_q(\nu)}{D^0_q(\L)}\label{xi*}
\eeq
where $m_\nu^\pm$ are the multiplicities of $V(\nu)$ in $W_\pm$, respectively,
so that
\beq
m_\nu=m_\nu^++m_\nu^-.
\eeq
\vskip.1in
\noindent {\bf Note:} in the case that $\L$ is typical it
necessarily follows that all
$V(\nu)$ in the tensor product decomposition (\ref{decomposition}) are
also typical so that (\ref{xi*}) is always well defined (c.f. lemma
\ref{lemma-unitary}).

\begin{Theorem}\label{link-th}
Consider the braid group $B_M$ and a braid $\t$ of the following
general form:
\beq
\t=\left (\s_{i_1}\right )^{k_1}\left (\s_{i_2}\right )^{k_2}\cdots
   \left (\s_{i_{M-1}}\right )^{k_{M-1}},~~~~~~~k_i\in{\bf Z}\label{theta}
\eeq
with $\{i_1,i_2,\cdots, i_{M-1}\}$ an arbitrary permutation of
$\{1,2,\cdots,M-1\}$. Then the following functional is a link polynomial
\beq
L(\hat{\t})=q^{-(\L,\L+2\rho)\sum_{i=1}^{M-1}k_i}\,
  \prod_{i=1}^{M-1}\xi^\L_{k_i}.\label{link-theta}
\eeq
In the case that $\L$ is typical,  $\xi_k^\L$ is given by (\ref{xi*}).
\end{Theorem}

\vskip.1in
\noindent {\bf VII. New Two-variable Link Polynomials}
\vskip.1in
We will now apply the technique developed in previous sections to
develop a general method for obtaining two-variable link polynomials
corresponding to any real $\L\in D^+$.
Again we restrict to the type-I quantum
superalgebras $\G=gl(m|n)$ or $\G=osp(2|2n)$.

Corresponding to any {\em real} $\L\in D^+$ we have the one-parameter family
of irreps
\bea
&&V(\L_\a)\equiv V(\L+\a\d),~~~~\a\in {\bf R},\no\\
&&\d=\left \{
\begin{array}{ll}
\sum_i\d_i, & ~~{\rm for}~ \G=gl(m|n)\\
\e_0, & ~~{\rm for}~ \G=osp(2|2n).
\end{array}
\right .
\eea
The module $V(\L_\a)$ is typical and unitary
for $|\a|$ sufficiently large. For example,
for the case $\G=gl(m|n)$, we have from proposition \ref{prop1}
a type (1) unitary module for
$\a> -(\L+\rho, \e_m-\d_n)=n-1-(\L,\e_m-\d_n)$,
and a type (2) unitary module for
$\a< -(\L+\rho, \e_1-\d_1)=1-m-(\L,\e_1-\d_1)$.
Below we asume $\a$ belongs to this range (although the final formula for
link polynomials should apply, by analytic continuation, to all real $\a$).

Here we obtain a representation of the braid generator $\s\in {\rm End}
[V(\L+\a\d)\otimes V(\L+\a\d)]$ and a formula for two variable link
polynomials. Consider the
$U_q(\G_0)$-module direct sum decomposition
\beq
V_0(\L)\otimes K(\L)=\bigoplus_\nu m_\nu V_0(\nu)\label{kac-decom}
\eeq
where $\G_0$ is the even subalgebra of $\G$ and $V_0(\L)$ the maximal
${\bf Z}$-graded component of $V(\L)$. Then for $|\a|$ sufficient large
(i.e. in the range considered above) we have the easily established
decomposition
\beq
V(\L+\a\d)\otimes V(\L+\a\d)=\bigoplus_\nu m_\nu V(\nu+2\a\d).
\eeq
Note that this decomposition may be obtained solely from a knowledge of
the $U_q(\G_0)$ modules occurring in $K(\L)$ and $U_q(\G_0)$ tensor
product rules. In principal this follows from the known characters of
$K(\L)$ and $V_0(\L)$.

 From our previous results we have the Casimir invariants
\beq
C^\L_k=(I\otimes {\rm str})[I\otimes \pi_{\L+\a\d}(q^{-2h_\rho})]\s^k
\eeq
which, from (\ref{xi*}), take the following eigenvalues
on $V(\nu+\a\d)$:
\beq
\xi^\L_k(q,\a)=q^{-kC(\L+\a\d)}\sum_\nu(-1)^{[\nu]}\,q^{\frac{k}{2}
  C(\nu+2\a\d)}\left (m_\nu^++(-1)^km_\nu^-\right )
  \frac{\chi_{\L+\a\d}(\Gamma_0)}{\chi_{\nu+2\a\d}
  (\Gamma_0)}\cdot \frac{D^0_q(\nu)}{D^0_q(\L)}\label{xi**}
\eeq
where use has been made of the fact that $\a\d$ is orthogonal to all even
roots and $\L+\a\d, ~\nu+2\a\d$ are all typical for $\a$ in the range
considered.

Now for $\t$ a braid of the general form (\ref{theta}), we arrive at
at the link polynomial
\beq
L(\hat{\t})=q^{-(\L+\a\d,\L+\a\d+2\rho)\sum_{i=1}^{M-1}k_i}\,
  \prod_{i=1}^{M-1}\xi^\L_{k_i}(q,\a)
\eeq
with $\xi^\L_k(q,\a)$ given by
(\ref{xi**}). In this way we obtain a two-variable
link polynomial corresponding to any real $\L\in D^+$.

\vskip.1in
\noindent {\bf Two-variable Link Polynomials from
  $U_q[gl(m|n)]$}
\vskip.1in
Following \cite{Del94b}, we assume $m\geq n$ and for $0\leq N\leq mn$
we call a Young diagram
$[\l]=[\l_1,\l_2,\cdots,\l_t],~\l_1\geq\l_2\cdots\geq\l_t\geq 0$ for
the permutation group $S_{N}$ (i.e. $\l_1+\l_2+\cdots+\l_t=N$)
%associated with $gl(m|n)$
{\bf allowable}, if it has at most $n$ columns and $m$ rows; i.e.
$t\leq m,~\l_i\leq n$.
%\vskip.1in
Associated with each such Young diagram $[\l]$
we define a weight of $gl(m|n)$
\beq
\L_{[\l]}=(\dot{0}_{m-t},-\l_t,\cdots,-\l_1|\underbrace{t,\cdots,t}_{\l_t},
  \underbrace{t-1,\cdots,t-1}_{\l_{t-1}-\l_t},\cdots,
  \underbrace{1,\cdots,1}_{\l_1-\l_2},
  \underbrace{0,\cdots,0}_{n-\l_1}).\label{weight}
\eeq
Using the basis $\{\e_i,~\d_j\}$, the weight $\L_{[\l]}$ may be expresed as
\beq
\L_{[\l]}=-\sum_{i=1}^t\l_i\e_{m-i+1}+t\sum_{j=1}^{\l_t}\d_j
  +\sum_{s=1}^t(t-s)\sum^{\l_{t-s}}_{j=\l_{t-s+1}}\d_j.\label{weight*}
\eeq

Let us consider the one-parameter family of finite-dimensional irreducible
$U_q(gl(m|n))$-modules $V(\L_\alpha)$ with highest weights of the form
$\Lambda_\alpha=(0,\cdots,0|
\alpha,\cdots,\alpha)\equiv (\dot{0}|\dot{\a})=\a\d$.
[That is the case $\L=(\dot{0}|\dot{0})$.]
These irreps $V(\a\d)$ are unitary of type (1) if $\a>n-1$ and
unitary of type (2) if $\a<1-m$. As mentioned above we assume real
$\a$ satisfying one of these conditions, in which case $V(\a\d)$ is
also typical of dimension $2^{mn}$.

We have the following decomposition of $V(\a\d)$
into irreps of the even subalgebra
$gl(m)\bigoplus gl(n)$:
\beq
V(\a\d)=\bigoplus^{mn}_{N=0}\bigoplus_{[\l]\in S_{N}}
  V_0(\L_{[\l]}+\a\d)
\eeq
where the summation is over allowed $N$-box Young
diagrams. Note that the index $N$ gives the ${\bf Z}$-graded level of
the irrep concerned. Alternatively we may simply write
\beq
V(\a\d)=\bigoplus_{[\l]}
  V_0(\L_{[\l]}+\a\d).
\eeq
The number of boxes $N_\l$ in the Young diagram $[\l]$ then gives the level.
We can deduce the tensor product decomposition
\beq
V(\a\d)\otimes V(\a\d)=\bigoplus_{[\l]}
  V(\L_{[\l]}+2\a\d).\label{glmn-decomposition}
\eeq
The parity of the module $V(\L_{[\l]}+2\a\d)$ is $(-1)^{N_\l}$.
The eigenvalue of the second order Casimir on the irrep $V(\L_{[\l]}+
2\a\d)$ can be shown to be
\bea
&&C(\L_{[\l]}+2\a\d)=2\sum_{i=1}^t\l_i(\l_i+1-2\a-2i)-2\a n(2\a+m),\nonumber\\
&&C(\a\d)=-\a n (\a+m).
\eea
Introduce the notation
\beq
\g_\a[\l]\equiv \frac{1}{2} C(\L_{[\l]}+2\a\d)-C(\a\d)=
  2\sum_{i=1}^t\l_i(\l_i+1-2\a-2i)-\a n(3\a+m).
\eeq
For $\t$ a braid of the general form (\ref{theta}) we arrive at the two
variable link polynomial
\beq
L(\hat{\t})=q^{-n\a(\a+m)\sum_{i=1}^{M-1}k_i}\,
   \prod_{i=1}^{M-1} \xi_{k_i}(q,\a)
\eeq
where now
\beq
\xi_k(q,\a)=\sum_{[\l]}(-1)^{(k-1)N_\l}\,q^{k\g_\a[\l]}
  \frac{\chi_{\a\d}(\Gamma_0)}{\chi_{\L_{[\l]}+2\a\d}(\Gamma_0)}\cdot
  \frac{D^0_q(\L_{[\l]}+2\a\d)}{D^0_q(\a\d)}.
\eeq
In this formula, the sum is again over all allowable Young diagrams.
This formula can be made fully explicit if we make use of
the easily established result (which takes a bit of algebra)
\bea
\chi_{\a\d}(\Gamma_0)\cdot \prod_{\b\in \Phi_1^+}[(\rho,\b)]_q
  &=&\prod_{i=1}^m\prod_{j=1}^n[i-j+\a]_q\nonumber\\
\chi_{\L_{[\l]}+2\a\d}(\Gamma_0) \cdot
  \prod_{\b\in \Phi_1^+}[(\rho,\b)]_q
  &=&\prod_{i=1}^m\prod_{j=1}^n[i-j-\l_i+2\a]_q
  \prod_{l=1}^t\frac{[\l_i-i-2\a+1-l]_q}{[\l_l+\l_i-i-2\a-l+1]_q}
\eea
where, in this last formula, it is implicitly understood that
$\l_i=0$ for $m\geq i>t$. We thus obtain
\beq
\xi_k(q,\a)=\sum_{[\l]}(-1)^{(k-1)N_\l}\,q^{k\g_\a([\l])}\,\chi_\a([\l])
  \,D^0_q(\L_{[\l]})
\eeq
where
\beq
\chi_\a([\l])\equiv \frac{\chi_{\a\d}(\Gamma_0)}{\chi_{\L_{[\l]}+2\a\d}
  (\Gamma_0)}=\prod_{i=1}^m\prod_{j=1}^n\frac{[i-j+\a]_q}{[i-j-\l_i+2\a]_q}
  \prod_{l=1}^t\frac{[\l_l+\l_i-i-2\a+1-l]_q}{[\l_i-i-2\a+1-l]_q}.
%&&D_q^0(\L_{[\l]})=\prod_{i<j}^m\frac{[\l_{m-j+1}\d_{j>m-t}-\l_{m-i+1}
%  \d_{i>m-t}+j-i]_q}{[j-i]_q}\prod_{k<l}^n\frac{[t(\d_{l<\l_t}-\d_{k<\l_t})
%  +\sum_{s=1}^t(t-s)(\d_{l<\l_{t-s}-\l_{t-s+1}}-\d_{k<\l_{t-s}-\l_{t-s+1}})
%  +k-l]_q}{[k-l]_q}
\eeq

\newsavebox{\yo}
\savebox{\yo}(4,8){\begin{picture}(8,8)(-2,0)
\put(0,0){\framebox(4,4){}}
\end{picture}}

\newsavebox{\yt}
\savebox{\yt}(4,8){\begin{picture}(8,11)(-2,0)
\multiput(0,0)(0,4){2}{\framebox(4,4){}}
\end{picture}}

\newsavebox{\yr}
\savebox{\yr}(4,10){\begin{picture}(8,12)(-2,3)
\multiput(0,0)(0,4){3}{\framebox(4,4){}}
\end{picture}}

\newsavebox{\yto}
\savebox{\yto}(10,8){\begin{picture}(10,11)(-2,0)
\multiput(0,0)(0,4){2}{\framebox(4,4){}}
\put(4,4){\framebox(4,4){}}
\end{picture}}

\newsavebox{\yro}
\savebox{\yro}(10,10){\begin{picture}(10,12)(-2,3)
\multiput(0,0)(0,4){3}{\framebox(4,4){}}
\put(4,8){\framebox(4,4){}}
\end{picture}}

\newsavebox{\yrt}
\savebox{\yrt}(10,10){\begin{picture}(10,12)(-2,3)
\multiput(0,0)(0,4){3}{\framebox(4,4){}}
\multiput(4,4)(0,4){2}{\framebox(4,4){}}
\end{picture}}

\newsavebox{\ytt}
\savebox{\ytt}(10,8){\begin{picture}(10,11)(-2,0)
\multiput(0,0)(0,4){2}{\framebox(4,4){}}
\multiput(4,0)(0,4){2}{\framebox(4,4){}}
\end{picture}}

\newsavebox{\yoo}
\savebox{\yoo}(10,8){\begin{picture}(10,8)(-2,0)
\put(0,0){\framebox(4,4){}}
\put(4,0){\framebox(4,4){}}
\end{picture}}

As an illustration, let us consider some specific cases in the remaining
part of this subsection.
\vskip.1in
\noindent{\bf Example (1): $U_q[gl(2|2)]$}
\vskip.1in

The tensor product decomposition is
\bea
V(\a\d)\otimes V(\a\d)&=&V(0,0|2\a, 2\a)\oplus
V({0},-1|2\a+1,2\a)\oplus
V(-1,-1|2\a+2,2\a)\nonumber\\
&&\oplus
V(0,-2|2\a+1,2\a+1)\oplus
V(-1,-2|2\a+2,2\a+1)\nonumber\\
&&\oplus
V(-2,-2|2\a+2,2\a+2).
\eea
We have in this case (using the Young diagram notation)
\bea
&&\g_\a(\cdot)=-2\a^2,~~~~
  \g_\a(\usebox{\yo})=-2\a(\a+1),\nonumber\\
&&\g_\a(\usebox{\yoo})=-2(\a^2+2\a-1),~~~~
  \g_\a(\usebox{\yt})=-2(\a+1)^2,\nonumber\\
&&\g_\a(\usebox{\yto})=-2\a(\a+3),~~~~
  \g_\a(\usebox{\ytt})=-2\a(\a+4),\nonumber\\
&&D^0_q(\cdot)=1,~~~~
  D^0_q(\usebox{\yo})=[2]_q^2,~~~~
  D^0_q(\usebox{\yoo})=[3]_q,\nonumber\\
&&D^0_q(\usebox{\yt})=[3]_q,~~~~
  D^0_q(\usebox{\yto})=[2]_q^2,~~~~
  D^0_q(\usebox{\ytt})=1,
\eea
while the $\chi_\a$ factors read
\bea
&&\chi_\a(\cdot)=\frac{[\a]_q^2[\a-1]_q[\a+1]_q}{[2\a]_q^2[2\a-1]_q[2\a+1]_q}
  ,\nonumber\\
&&\chi_\a(\usebox{\yo})=\chi_\a(\usebox{\yto})=
  \frac{[\a]_q^2[\a-1]_q[\a+1]_q}{[2\a]_q^2[2\a-2]_q[2\a+2]_q},\nonumber\\
&&\chi_\a(\usebox{\yoo})=
  \chi_\a(\usebox{\yt})=
  \frac{[\a]_q^2[\a-1]_q[\a+1]_q}{[2\a+1]_q[2\a-1]_q[2\a-2]_q[2\a+2]_q}
  ,\nonumber\\
&&\chi_\a(\usebox{\ytt})=\chi_\a(\cdot)=
  \frac{[\a]_q^2[\a-1]_q[\a+1]_q}{[2\a]_q^2[2\a-1]_q[2\a+1]_q}.
\eea
It follows that
\bea
\xi_k(q,\a)&=&q^{k\g_\a(\cdot)}\chi_\a(\cdot)D^0_q(\cdot)+
  q^{k\g_\a(\usebox{\ytt})}\chi_\a(\usebox{\ytt})D^0_q(\usebox{\ytt})
  \nonumber\\
& &-(-1)^k\left [q^{k\g_\a(\usebox{\yo})}\chi_\a(\usebox{\yo})D^0_q(\usebox
  {\yo})+
  q^{k\g_\a(\usebox{\yto})}\chi_\a(\usebox{\yto})D^0_q(\usebox{\yto})\right]
  \nonumber\\
& &+q^{k\g_\a(\usebox{\yoo})}\chi_\a(\usebox{\yoo})D^0_q(\usebox
  {\yoo})+
  q^{k\g_\a(\usebox{\yt})}\chi_\a(\usebox{\yt})D^0_q(\usebox{\yt})
  \nonumber\\
&=&q^{-2k\a(\a+2)}\frac{(q^{4k\a}+q^{-4k\a})[\a+1]_q[\a-1]_q}{(q^\a+q^{-\a})^2
    [2\a+1]_q[2\a-1]_q}\nonumber\\
& &-(-1)^k q^{-2k\a(\a+2)}\frac{(q^{2k\a}+q^{-2k\a})
    [2]_q^2}{(q^\a+q^{-\a})^2(q^{\a-1}+q^{-\a+1})(q^{\a+1}+q^{-\a-1})}
    \nonumber\\
& &+\frac{q^{-2k(\a^2+\a+1)}\,(q^{2k\a}+q^{-2k\a})[3]_q
   (q^\a+q^{-\a})^2}{(q^{2\a-1}-q^{-2\a+1})(q^{2\a+1}-q^{-2\a-1})
   (q^{\a-1}+q^{-\a+1})(q^{\a+1}+q^{-\a-1})}.
\eea

\vskip.1in
\noindent{\bf Example (2):} $U_q[gl(m|1)]$
\vskip.1in

We have the tensor product decomposition
\bea
V(\a\d)\otimes V(\a\d)&=&V(\dot{0}|2\a)\oplus
V(\dot{0},-1|2\a+1)\no\\
&&\oplus V(\dot{0},-1,-1|2\a+2)\oplus \cdots
\oplus V(-\dot{1}|2\a+m).
\eea
In this case $D^0_q(\L_{[\l]})$ reads
\beq
D^0_q(\L_{[\l]})=\prod_{i=1}^t\frac{[m+1-i]_q}{[t+1-i]_q}\equiv
  \frac{[m]_q!}{[m-t]_q![t]_q!}
\eeq
and $\g_\a[\l],~~\chi_\a([\l])$ reduce to, respectively,
\bea
&&\g_\a[\l]=-t(t-1)-\a(\a+2t),\no\\
&&\chi_\a([\l])=\prod_{i=1}^m\frac{[i+\a-1]_q}{[i+2\a+t-1-\l_i]_q}.
\eea
The $\xi_k(q,\a)$ have the following form,
\beq
\xi_k(q,\a)=\sum_{t=0}^m(-1)^{(k-1)t}\,q^{-k[t(t-1)+\a(\a+2t)]}\,
  \prod_{i=1}^t\frac{[m+1-i]_q[i+\a-1]_q}{[t+1-i]_q[i+2\a+t-2]_q}
  \prod_{i>t}^m\frac{[i+\a-1]_q}{[i+2\a+t-1]_q}.
\eeq

\vskip.1in
\noindent {\bf Two-variable Link Polynomials from Adjoint Representation
   of $U_q[gl(2|1)]$}
\vskip.1in
As another illustration of how the general formalism works it is instructive to
consider the case $\L=\psi,~~\psi=(1,0|-1)$ the highest weight of the
adjoint representation of $gl(2|1)$.
This example is of interest since it affords the
simplest example of a two-variable link polynomial in which a multiplicity
occurs in the tensor product space.

First note that in this case $\e_1-\e_2$ is the single even positive root
and $\e_1-\d_1,~\e_2-\d_1$ are the two odd positive roots, from which we deduce
that for any $\L=(\L_1,\L_2|\bar{\L}_1)$
\beq
D^0_q[\L]=[\L_1-\L_2+1]_q,~~~~~\chi_\L(\Gamma_0)=[\L_1+\bar{\L}_1+1]_q
  [\L_2+\bar{\L}_1]_q.
\eeq

For the Kac-module $K(\psi)$ we have the $U_q(\G_0)$-module
($\G_0=gl(2)\bigoplus u(1)$) decomposition (illustrated in terms of
${\bf Z}$-graded levels):
\beq
K(\psi)=V_0(1,0|-1)\bigoplus V_0(1,-1|0)\bigoplus V_0(0,0|0)
  \bigoplus V_0(0,-1|1)
\eeq
which is easily seen to be $2^2\cdot 2=8$ dimensional as required. Thus
\bea
V_0(\psi)\otimes K(\psi)&=&V_0(1,0|-1)\otimes V_0(1,0|-1)\nonumber\\
& &\bigoplus V_0(1,0|-1)\otimes\left [V_0(1,-1|0)\oplus V_0(0,0|0)\right]\no\\
& &\bigoplus V_0(1,0|-1)\otimes V(0,-1|1)\no\\
&=&V_0(2,0|-2)\bigoplus V_0(1,1|-2)\bigoplus V_0(2,-1|-1)\no\\
& &\bigoplus 2\,V_0(1,0|-1)\bigoplus V_0(1,-1|0)\bigoplus V_0(0,0|0)
\eea
which yields the tensor product decomposition:
\bea
V(\psi+\a\d)\otimes V(\psi+\a\d)&=&V(2,0|2\a-2)\bigoplus
  V(1,1|2\a-2)\bigoplus V(2,-1|2\a-1)\no\\
& &\bigoplus 2\,V(1,0|2\a-1)\bigoplus V(1,-1|2\a)\bigoplus V(0,0|2\a).
\eea
It is seen that $V(1,0|2\a-1)$ occurs twice
in the tensor product space. From the above
${\bf Z}$ gradation on $V_0(\psi)\otimes K(\psi)$ we obtain
\beq
(-1)^{[\nu]}=\left \{
\begin{array}{ll}
-1, & {\rm for}~ \nu=(2,-1|2\a-1),~~(1,0|2\a-1)\\
1, & {\rm otherwise}.
\end{array}
\right .
\eeq
In the $q\rightarrow 1$ limit the above tensor product module decomposes
into symmetric and anti-symmetric components (which determine the
parities):
\beq
V(\psi+\a\d)\otimes V(\psi+\a\d)=W_+\bigoplus W_-
\eeq
with
\bea
W_-&=&V(1,1|2\a-2)\bigoplus V(2,-1|2\a-1)\bigoplus V(1,0|2\a-1)
   \bigoplus V(0,0|2\a),\no\\
W_+&=&V(2,0|2\a-2)\bigoplus
  V(1,0|2\a-1)\bigoplus V(1,-1|2\a).
\eea
Note that there is one copy of $V(1,0|2\a-1)$ in each of these spaces.
For the Casimirs we have
\beq
\frac{1}{2}C(\nu+2\a\d)-C(\psi+\a\d)=\left\{
\begin{array}{ll}
-\a(\a+2), & \nu=(2,0|-2), ~~(1,-1|0)\\
-(\a^2+1), & \nu=(1,0|-1)\\
-(\a^2+2\a+2), & \nu=(0,0|0)\\
-(\a^2-2\a+2), & \nu=(1,1|-2)\\
-\a^2+2, & \nu=(2,-1|-1).
\end{array}
\right .
\eeq

Collecting together all of this information and substituting into
(\ref{xi**}) we arrive at
\bea
\xi^\psi_k(q,\a)&=&q^{-k\a(\a+2)}\frac{[\a+1]_q[\a-1]_q[3]_q}{[2\a+1]_q
   [2\a-2]_q[2]_q}+q^{-k\a(\a+2)}
   \frac{[\a+1]_q[\a-1]_q[3]_q}{[2\a+2]_q[2\a-1]_q[2]_q}\no\\
& & +(-1)^kq^{-k(\a^2+2\a+2)}\frac{[\a+1]_q[\a-1]_q}{[2\a+1]_q[2\a]_q[2]_q}
   +(-1)^kq^{-k(\a^2-2\a+2)}\frac{[\a+1]_q[\a-1]_q}{[2\a]_q[2\a-1]_q[2]_q}\no\\
& &-(-1)^kq^{-k(\a^2-2)}\frac{[\a+1]_q[\a-1]_q[4]_q}
   {[2\a+2]_q[2\a-2]_q[2]_q}\no\\
& &-\left (1+(-1)^k\right )q^{-k(\a^2+1)}\frac{[\a+1]_q[\a-1]_q}
   {[2\a+1]_q[2\a-1]_q}.
\eea

\vskip.1in
\noindent {\bf Two-variable Link Polynomials from $U_q[osp(2|2n)]$}
\vskip.1in
Consider the one-parameter family of $2^{2n}$-dimensional irreducible
$U_q[osp(2|2n)]$-modules $V(\L_\a)$ with highest weights of form
$\L_\a=(\a|0,\cdots,0)\equiv \a\e_0$ (and with lowest weight $\L^-_\a=
(\a-2n)\e_0$).
$V(\a\e_0)$ is unitary and typical  provided
that $\a<0$ or $\a>2n$. We therefore consider the tensor product module
$V(\a\e_0)\otimes V(\a\e_0)$ with $\a<0$ or $\a>2n$ which decomposes as
\beq\label{label}
V(\a\e_0)\otimes V(\a\e_0)=\bigoplus_{c=0}^n\bigoplus_{d=0}^{n-c}
   V(\L_{c,d})\label{osp-decom}
\eeq
with
\beq
\L_{c,d}=(2\a-c-2d)\e_0+\l_c,~~~~~ \l_c=
\sum_{i=1}^c\e_i.
\eeq
The decomposition (\ref{osp-decom}) is obtained from known character
formulae [c.f. eq.(\ref{kac-decom})].

{}From the ${\bf Z}$ gradation on $V(\a\e_0)$ we can deduce that the level
of the module $V(\L_{c,d})$ is equal to $c+2d$.
Thus the parity of the module $V(\L_{c,d})$ is
$(-1)^{c+2d}$. The Casimir eigenvalues read
\bea
&&C(\L_{c,d})=4(\a-d)(n+c+d-\a)-2c(c-1)\no\\
&&C(\a\e_0)=\a(2n-\a).
\eea
For $\t$ a braid of the general form (\ref{theta}) we thus arrive at the two
variable link polynomial
\beq
L(\hat{\t})=q^{-\a(2n-\a)\sum_{i=1}^{M-1}k_i}\,
  \prod_{i=1}^{M-1} \xi_{k_i}(q,\a)
\eeq
where the $\xi_k(q,\a)$'s are given by
\beq
\xi_k(q,\a)=\sum_{c=0}^n\sum_{d=0}^{n-c}(-1)^{(k-1)(c+2d)}\,q^{k\g_\a}
  \frac{\chi_{\a\e}(\Gamma_0)}{\chi_{\L_{c,d}}(\Gamma_0)}\cdot
  \frac{D^0_q(\L_{c,d})}{D^0_q(\a\e)}
\eeq
with
\beq
\g_\a\equiv \frac{1}{2}C(\L_{c,d})-C(\a\e).
\eeq
%In view of the fact that $\Phi_0^+=\{\d_i-\d_j,~\d_i+\d_j,~2\d_i\;|\; i<j\},~
%\Phi_1^+=\{\e\pm\d_i\}$,
After a bit algebra, we end up with
\bea
\chi_{\a\e}(\Gamma_0) \cdot \prod_{\b\in\Phi_1^+}[(\rho,\b)]_q&=&
  \prod_{i=1}^n [2n+1-i-\a]_q[i-\a-1]_q,\no\\
\chi_{\L_{c,d}}(\Gamma_0) \cdot
  \prod_{\b\in\Phi_1^+}[(\rho,\b)]_q
&=&\prod_{i=1}^n [c+2d+2n+1-i-2\a-\d_{i\leq c}]_q\no\\
& &\cdot [c+2d+i-2\a-1+\d_{i\leq c}]_q
\eea
where $\d_{i\leq c}$ equals 1 for $i\leq c$ and zero otherwise.
We thus obtain
\beq
\xi_k(q,\a)=\sum_{c=0}^n\sum_{d=0}^{n-c}(-1)^{(k-1)(c+2d)}\,q^{k\g_\a}
  \chi_\a(c,d) \cdot D^0_q(\l_c)
\eeq
where
\bea
\chi_{\a}(c,d)\equiv \frac{\chi_{\a\d}(\Gamma_0)}
  {\chi_{\L_{c,d}}(\Gamma_0)}&=&
  \prod_{i=1}^n \frac{[2n+1-i-\a]_q[i-\a-1]_q}
  {[c+2d+2n+1-i-2\a-\d_{i\leq c}]_q
  [c+2d+i-2\a-1+\d_{i\leq c}]_q},\no\\
D^0_q(\l_c)&=&\prod_{i<j}^c\frac{[2(n+2)-i-j]_q}{[2(n+1)-i-j]_q}
  \,\prod_{l=1}^c\frac{[2(n+2-l)]_q}
  {[2(n+1-l)]_q}.
\eea

\vskip.1in
\noindent {\bf VIII. Discussion}
\vskip.1in
We have demonstrated how link polynomials can be constructed associated with
any finite dimensional unitary irrep of a type-I quantum superalgebra.
This is achieved by successfully overcoming a fundamental problem in
computing the eigenvalues of Casimir invariants for the quantum superalgebras.
Applying our results to one-parameter families of inequivalent irreps,
we have been able to construct infinite families of non-equivalent
two-variable link polynomials. Such two-variable link polynomials were
previously known only for some isolated cases. For a class of braids,
we have computed the link polynomials in fully explicit form.

\vskip.3in
\begin{center}
{\bf Acknowledgements:}
\end{center}
The financial support from the Australian Research Council is gratefully
acknowledged. Y.-Z.Z is financially supported by the Kyoto University
Foundation.

%\newpage

\end{document}